\documentclass[11  pt]{article}
\usepackage{amsmath,amsfonts,amssymb,amsthm}

\usepackage[dvips]{graphics}
\usepackage[english]{babel}
\usepackage{makeidx}
\usepackage{latexsym}
\usepackage[latin1]{inputenc}
\newtheorem{theo}{Theorem}[section]
\newtheorem{lema}{Lemma}[section]
\newtheorem{prop}{Proposition}[section]
\newtheorem{defi}{Definition}[section]
\newtheorem{example}{Example}[section]
\newtheorem{rema}{Remark}[section]
\newtheorem{cor}{Corollary}[section]

\begin{document}
 \baselineskip=1.0\baselineskip
  \title{Bj\"orling problem for timelike surfaces in the Lorentz-Minkowski space}
    \author{R. M. B. Chaves $^*$, \,\,M. P. Dussan
\footnote{Universidade de S\~ao Paulo, Departamento de
Matem\'atica -IME. CEP: 05508-090. S\~ao Paulo. Brasil}\,, M.
Magid \footnote{Department of Mathematics, Wellesley College,
Wellesley, MA 02458}}
\date{}
  \maketitle

 \begin{abstract}
We introduce a new approach to the study of  timelike minimal surfaces in the Lorentz-Minkowski space through a split-complex representation formula for this kind of surface. As applications, we solve the Bj\"orling problem for timelike surfaces and obtain interesting examples and related results. Using the Bj\"orling representation, we also obtain characterizations of minimal timelike surfaces of revolution as well as of minimal ruled timelike surfaces in the Lorentz-Minkowsi space.   

\end{abstract}
 MSC 53A10
 \section{Introduction}

 Minimal surfaces play an important role in Differential Geometry 
and  also in Physics, especially in  problems related to  General Relativity. In Euclidean space, the
classical Björling problem (see \cite{G}, \cite{DHKW}) was proposed by Björling \cite{B} in 1844 and consists of the construction of a minimal surface in $\mathbb R^3$ containing
a prescribed analytic strip. The solution to this problem was obtained by
Schwarz in \cite{S}. These results have inspired many authors to obtain analogous results in other ambient spaces. For instance, in \cite{ACG},  J. Aledo, R. Chaves and A. Gálvez studied the Bj\"orling problem  in the context of affine geometry, and  more recently Aledo, Martinez and Milan in \cite{AMM} generalized the results in  \cite{ACG} when solving the Bj\"orling problem for affine maximal surfaces. In \cite{AV}, Asperti and Vilhena studied the problem for spacelike surfaces in $\mathbb L^4$. Moreover,  F. Mercuri and I. Onnis, in \cite{MO}, studied the problem  when the ambient space is a Lie group, while Aledo, Galvez and Mira, in \cite{AGM}, solved a Bjorling-type problem for flat surfaces in the 3-sphere. For other references see \cite{GM} and \cite{KY}. 

 The local geometry of surfaces in the Lorentz-Minkowski space $\mathbb L^3$ is more complicated than that of surfaces in $\mathbb R^3$, since in $\mathbb L^3$ the  vectors have different causal characters which yield more cases to be analysed. Hence we could consider spacelike, timelike or lightlike surfaces in $\mathbb L^3$.  In this paper we consider timelike minimal surfaces in  $\mathbb L^3$.  Although timelike minimal surfaces neither maximize nor minimize surface area, they have nice geometric properties similar to minimal surfaces in Euclidean space $\mathbb R^3$. For instance, they also admit an Enneper-Weierstrass representation, introduced by M. Magid in \cite{Mag5}. But it is well known there are many differences in their behaviour.

 The Bj\"orling problem for spacelike surfaces in  $\mathbb L^3$ was considered by L. Alías, R. Chaves and P. Mira in  \cite{ACM} via a complex representation formula, while in \cite{MP} Mira and Pastor proposed  the problem of establishing a Bj\"orling type formula for timelike minimal surfaces. As it is well known and as was pointed in \cite{MP}, any general method for studying timelike minimal surfaces seems to be of interest, since this theory is not much developed. So, we were motivated to investigate  the Bj\"orling  problem for the timelike surfaces in $\mathbb L^3$. Our approach considers  the ring of split-complex numbers, denoted henceforth by $\mathbb C'$, which plays a role similar to that played by the ordinary complex-numbers in the spacelike case.  We note that the split-complex analysis and split-complex geometry depending on the split-complex numbers have been appropriately used to study  timelike surfaces, as well as their applications in physics; see for instance \cite{Er}, \cite{FI} and \cite{IT}. 
 
In this paper we  assume that all the real-valued  functions $\gamma(t)$ are real analytic. This condition allows us, after extending the function $\gamma$ to a subset $\mathcal O \subset \mathbb C'$,  to construct appropiate split-holomorphic functions defined in $\mathcal O$, whose real part will represent an analytic solution of the Bj\"orling problem in the timelike setting.  In others words we construct split-complex representation formulas for analytic timelike minimal surfaces that are solutions of the Bj\"orling problem. Since in the timelike setting we need to consider the causal character of the analytic curves, there are two Bj\"orling problems which we have to  study: Assuming that the analytic strip contains either an analytic timelike curve or an analytic spacelike curve  $\gamma: I \to \mathbb L^3$, and that an orthogonal unit analytic spacelike vector field $W$  is defined along $\gamma$.   We call these respectively the timelike or spacelike Bj\"orling problem. Then using our split-complex representation formulas we prove the existence and uniqueness of analytic solutions  to the two Bj\"orling problems as well as   obtaining important results and many interesting examples of minimal  surfaces with prescribed geometric properties. We note that since there is non-uniqueness of the solution of the Bj\"orling problem when considering the initial data $\gamma(t)$ as a null curve (as we can see from Example \ref{Exem} below),  we have just to consider the two cases above. 

\vspace{0.2cm}
We prove the following split-complex representation formula in order to solve the timelike Bj\"orling problem. A similar formula for the spacelike Bj\"orling problem is proven in Theorem \ref{eq:TB2}.

\begin{theo}(Timelike Bj\"orling representation)\label{eq:TB1}
Let $X:\frak O \subset \mathbb C'\to \mathbb L^3$ be a timelike
minimal surface in $\mathbb L^3$ and set $\gamma(t)=X(t,0)$, $W(t)
= N(t,0)$ on a real interval $I$ in $\frak O.$  Choose any simply
connected region $\frak R\subset \frak O$ containing $I$ over
which we can define split-holomorphic extensions $\gamma(z), W(z)$
for all $z\in \frak R$.  Then for all $z \in \frak R$ we have:
\begin{eqnarray}\label{eq:Bjorling}
X(z)=Re\left(\gamma(z) + k'\int_{t_o}^z W(w) \times \gamma'(w) dw
\right),
\end{eqnarray}
where $t_o$ is an arbitrary fixed point in $I$ and the path
integral is taken over any path in $\frak R$ from $t_o$ to $z$.
  \end{theo}
  
  Using the Bj\"orling approach  we also give alternative proofs of the well known characterizations of  minimal timelike surfaces of revolution and minimal timelike ruled surfaces in $\mathbb L^3$  (Woestijne \cite{W}).  More specifically, we prove the following results.
  
\begin{theo}\label{eq:rev} Every analytic minimal timelike surface of revolution in $\mathbb L^3$ is congruent to a part of one of the following surfaces:
 
 i) a Lorentzian elliptic catenoid.
 
 ii) a Lorentzian hyperbolic catenoid

 iii) a Lorentzian surface with spacelike profile curve.
 
 iv) a Lorentzian parabolic catenoid.
 
 \end{theo}
 
\begin{theo}\label{eq:ruled3}
Every analytic  minimal timelike ruled surface of $\mathbb L^3$ is congruent to a part of one of the following surfaces:

i) A Lorentzian plane

ii) The helicoid of the 1st kind.

iii)  The helicoid of the 2nd kind.

iv) The helicoid of the 3rd kind.

v) The conjugate surface of Enneper of the 2nd kind. 

vi) A flat $B$-scroll over a null curve.
\end{theo}

We have organized the paper  as follows. In Section 2 we fix the notation and give some preliminaries involving the split-complex numbers. Section 3 contains the key results regarding the analysis of extending real analytic functions to split-holomorphic functions.  In Section 3 we also state and solve the Bj\"orling problems determining minimal timelike surfaces in $\mathbb L^3$, constructed in terms of its Bj\"orling data.  Sections 4 and 5 contain, respectively, the characterization of analytic minimal timelike surfaces of revolution and analytic minimal timelike ruled surfaces, Theorems \ref{eq:rev} and \ref{eq:ruled3},  together with some  examples.

 \section{Preliminaries}
  
 Following the notation in \cite{IT}, we begin with a definition:

 \begin{defi} The split-complex numbers $\mathbb C' =\{t + k's | t, s \in \mathbb R, k'^2=1, 1k'=k'1\}$ is a commutative algebra over $\mathbb R$.  If $z=t + k's$ then $Re(z)=t$, $Im(z)=s$, $\bar{z}=t-k's$.  The indefinite metric on $\mathbb C'$ is given by $-z\bar{z}=-t^2 + s^2$.  \end{defi}
 \begin{defi} A function $f:\mathbb C' \to \mathbb C'$, $f(z)=f(t+k's)=u(t,s)+k'v(t,s)$ is split-holomorphic if and only if $u_t=v_s$ and $u_s=v_t$.  (See \cite{Er}, \cite{IT}). \end{defi}
Note that in this case $$\frac{d f}{d z}=f'(z) =\frac{1}{2}(\partial_t + k'
\partial_s)(u + k'v) = u_t + k' v_t = v_s + k' u_s.$$
 \begin{prop} If  $C$ is a curve in the $\mathbb C'$ plane and $f(z)$ is a split-holomorphic function on $C$ with a continuous derivative $f'(z)$, then
\begin{eqnarray}
\int_C f'(z)dz = f(z)|_C\end{eqnarray} and the integral is clearly
path independent.\end{prop}
\begin{proof}
We use the standard definition of a line integral.
\begin{eqnarray*} \int_C f'(z)dz =\int_C f'(t+k's)(dt + k' ds)=\\\int_Cu_t dt + u_s ds + k'(v_tdt + v_sds)
= \int_C d(u+k'v).
\end{eqnarray*}
\end{proof}
\begin{prop}
If $f=u + k'v $ is split holomorphic, then there is a
split-holomorphic function $g$ so that $g'=f.$  \end{prop}
\begin{proof}
We know that $u_t=v_s$ and $u_s=v_t$.  Take $ \beta,
\alpha:\mathbb R^2 \to \mathbb R$ so  that $\beta_s=u, \beta_t=v$
and $\alpha_t=u, \alpha_s=v$.  Let $g=\alpha + k' \beta$.  Then
$g'=\alpha_t + k' \beta_t= u + k' v= f.$
\end{proof}
Thus every line integral of a split-holomorphic function is path
independent.
\begin{prop}
If $f=u+k'v$ is a split-holomorphic function with $u, v \in C^2$
then $f'$ is again split-holomorphic.
\end{prop}
\begin{proof}
$f'(t,s)=u_t + k'v_t = v_s + k'u_s$.  We must show that
$u_{tt}=v_{ts}, u_{ts}=v_{tt}.$  We know that $u_t=v_s$ and
$u_s=v_t$, so we are done, as long as the mixed partials are equal.
\end{proof}

We  also note that using the following theorem and the Real Analytic Inverse Function Theorem (\cite{KP}), we can see that a split-holomorphic analytic function has a split-holomorphic analytic inverse if the determinant of the Jacobian is non-zero.

\begin{theo}  Let $U \subset \mathbb R^2$ be an open subset and let $f_1:U\to \mathbb R$ and $f_2 :U\to \mathbb R$ have continuous partial derivatives.  Consider the equations:
\begin{eqnarray}\label{eq:function}
f_1(t,s)=u\\
f_2(t,s)=v
\end{eqnarray}
near a given solution $(t_o,s_o)$, $(u_o,v_o)$.  If the Jacobian matrix $J_f(t_o,s_o)$ is non-zero then the equation (\ref{eq:function}) can be solved uniquely as $(t,s)=g(u,v)$ for $(t,s)$ near $(t_o,s_o)$ and $(u,v)$ near $(u_o,v_o)$.   Moreover, the function $g$ has continuous partial derivatives.
\end{theo}

\begin{defi}$\mathbb L^3$ is $\mathbb R^3$ with the indefinite inner product $$\left<(x_1,x_2,x_3), (y_1,y_2,y_3)\right> = -x_1y_1 + x_2y_2 + x_3y_3.$$ \end{defi}
Let $X:M^2_1\to \mathbb L^3$  be a timelike surface, i.e.,  a
surface which inherits a non-degenerate metric $h$ of signature
$(1,1)$ on every tangent space.  Following \cite{Wei}  we call the
pair $(M^2_1, [h])$ a Lorentz surface defined by $h$ where $[h]$ denotes the class of metrics conformal
to $h$ by a positive factor.
This is the analog of a Riemann surface in the timelike setting.

Let $z=t + k's$, where $t$ and $s$ are conformal coordinates and
$$\phi_j =\frac{\partial X_j}{\partial z} =
\frac{1}{2}\left(\frac{\partial X_j}{\partial t} +
k'\frac{\partial X_j}{\partial s}\right),$$ 
where $X_j$ represents a component of timelike immersion.

Observe that $$-\phi_1^2
+ \phi_2^2 + \phi_3^2 = \left <X_s,X_s\right> + \left<X_t, X_t\right> + 2 \left<X_t,
X_s\right> =0.$$  If we set $|a + k'b|^2 = b^2 - a^2$ then
$$-|\phi_1|^2 + |\phi_2|^2 + |\phi_3|^2 = \left<X_s, X_s\right> - \left<X_t, X_t\right>>0.$$

Consider the split-complex 1-forms defined by $\Phi_j = \phi_j
dz$. By looking at a conformal change of coordinates and using
Proposition 2.1 in \cite{Mag5}, we can see this defines a global
form on $M^2_1$.

In her unpublished thesis  \cite{Be} Berard shows that $X$ is
minimal if and only if $X_{jss}=X_{jtt}$ with respect to isothermal
coordinates $\{t,s\}$.  Thus the $\phi_j$ are split-holomorphic
if and only if $X$ is a minimal immersion in $\mathbb L^3$.

Note that if $\Phi_k$ are globally defined, then, following
\cite{FT}, pp. 77--78, we have, in a local coordinate patch:
\begin{eqnarray*}
2 Re \int_\gamma \phi_j dz = Re \int_\gamma
\frac{1}{2}\left(\frac{\partial X_j}{\partial t} +
k'\frac{\partial X_j}{\partial s}\right)(dt + k'ds)\\ =
\int_\gamma \left(\frac{\partial X_j}{\partial t}dt +
\frac{\partial X_j}{\partial s} ds\right)= \int_\gamma dX_j =
X_j\left |_\gamma \right..
\end{eqnarray*}
Thus the integral over any closed curve has real part zero.  The
converse is also true.
\begin{theo} \label{eq:immersion} Let $\Sigma$ be a Lorentzian surface and choose three split-holomorphic one-forms $\Phi_1$, $\Phi_2$, $\Phi_3$ globally defined on $\Sigma$ satisfying:
\begin{enumerate}
\item$-\Phi_1^2 + \Phi_2^2 + \Phi_3^2 =0$.
\item$-|\Phi_1|^2 + |\Phi_2|^2 + |\Phi_3|^2>0$.
\item Each $\Phi_j$ has no real periods.
\end{enumerate}
Then the map $X: \Sigma \to \mathbb L^3$ given by $$X(z)=2 Re
\int_{\gamma_z}(\Phi_1,\Phi_2,\Phi_3) \, dz,$$ where $\gamma_z$ is a
path from the fixed basepoint $z$ is a minimal immersion in
$\mathbb L^3$.\end{theo}

\begin{rema} We could also use the split-complex variable
$w=k'z=s+k't$ in the above formulas, setting
 $$\psi_j =\frac{\partial X_j}{\partial w} = \frac{1}{2}\left(\frac{\partial X_j}{\partial s} + k'\frac{\partial X_j}{\partial t}\right).$$  After replacing $\Phi_j$ by $\Psi_j$ the formulas are the same, except that $$-|\psi_1|^2 + |\psi_2|^2 + |\psi_3|^2 =  - \left< X_s, X_s \right> + \left< X_t, X_t \right> <0.$$ We will use this alternative choice of variable in Subsection 3.1.
\end{rema}

\section{Proofs of the main results}
Throughout  these last sections we assume that all the real-valued functions $\gamma(t)$ 
are $C^\omega$, i.e., they are real analytic (have power series
representations).  
\vspace{0.3cm}

 We begin by proving results about extending $\gamma(t)$ as a split-holomorphic function and  proving that a split-complex analytic function $f(z)$ on a domain $\mathcal D \subset   \mathbb C'$ is uniquely determined in $\mathcal D$ by knowledge of the derivatives $f^{(k)} (\alpha)$ at a point $\alpha \in \mathcal D$. 

\begin{prop}Let $\gamma(x)$ be a real analytic function given by $\gamma(x)=\sum_{k=0} ^ \infty \gamma_k x^k$ which converges in $|x| < R.$  This function can be extended split-holomorphically in a neighborhood of $0\in \mathbb C'$.  \end{prop}

\begin{proof}
First consider the new series
\begin{eqnarray}\label{eq:newseries}
\sum_{k=0}^\infty 2^{k-1}\gamma_ky^k.\end{eqnarray}   Using Hadamard's formula for the radius of convergence of a power series, we can see that the radius of convergence is $R/2$ (\cite{KP}).  We are considering 
\begin{eqnarray*} 
\gamma(t + k's) &=& u(t,s)+k'v(t,s) = \sum_{k=0}^\infty \gamma_k (t+k's)^k 
= \gamma_0 + \gamma_1 t + \gamma_2(t^2 + s^2)\\ &+&\gamma_3(t^3 + 3 ts^2) + \dots 
+k'(\gamma_1 s + \gamma_2(2ts) + \gamma_3(3t^2 s+ s^3)
+\dots.\end{eqnarray*} 
It is clear that both the real and
imaginary part converges on $|t| + |s|<R/2$, using the series
defined in (\ref{eq:newseries}).

Now we prove that $u_t=v_s$ and $u_s=v_t$. In fact, we just prove first equality on the even terms of $u(t,s)$ since  the other the equalities follow in a similar way. 

For $ k$ even, the $k^{th}$ term of $u(t,s)$ is $\gamma_k(t^k +
\begin{pmatrix} k\\2 \end{pmatrix}t^{k-2}s^2 + \begin{pmatrix}
k\\4 \end{pmatrix}t^{k-4}s^4+\dots + s^k).$ Then in  $u_t$ this
yields:
$$\gamma_k(kt^{k-1} +(k-2) \begin{pmatrix} k\\2 \end{pmatrix}t^{k-3}s^2 + (k-4)\begin{pmatrix} k\\4 \end{pmatrix}t^{k-5}s^4+\dots + 2 \begin{pmatrix} k\\k-2\end{pmatrix}ts^{k-2} ).$$

The $k^{th}$ term of $v(t,s)$ is:
 $\gamma_k(kt^{k-1}s + \begin{pmatrix} k\\3 \end{pmatrix}t^{k-3}s^3 + \begin{pmatrix} k\\5 \end{pmatrix}t^{k-5}s^5+\dots +  \begin{pmatrix} k\\k-1 \end{pmatrix}ts^{k-1}).$

 We can then see that $v_s$ has the  term:
  $$\gamma_k(kt^{k-1} + 3\begin{pmatrix} k\\3 \end{pmatrix}t^{k-3}s^2 + 5\begin{pmatrix} k\\5 \end{pmatrix}t^{k-5}s^4+\dots +  (k-1)\begin{pmatrix} k\\k-1 \end{pmatrix}ts^{k-2}).$$
  The proof finishes when we note that $$(k-2j)\begin{pmatrix} k\\2j \end{pmatrix} = (2j +1)\begin{pmatrix} k\\ 2j+1 \end{pmatrix}.$$
\end{proof}
\begin{defi}  Let $U\subset \mathbb R^2$ be an open subset. The function $f: U\subset \mathbb R^2 \to \mathbb R$ is called {\it{real analytic}} on $U$ if, for each $p\in U$, $f$ may be represented by a convergent power series in some neighborhood of $p$. We write \lq \lq$f$ is $C^\omega$".
\end{defi}
\begin{prop}\label{eq:3.2}
 If $\gamma(t)$ is a real-valued analytic function with two split-holomorphic extensions $u + k' v$ and $a + k' b$ satisfying  $a(t,0)=u(t,0)=\gamma(t)$ in an open set, then they agree everywhere they are defined.
 \end{prop}
\begin{proof} We can see, using analytic continuation of real analytic functions of one variable, (p. 14 (\cite{KP})), that $a(t,0)=u(t,0)=\gamma(t)$ for all $t$.  We know the extension of $\gamma(t)$, $f(t,s)=u(t,s) + k' v(t,s)$ is split-complex holomorphic.  Thus, following \cite{Er} if we let $$t=\frac{x+y}{2} \hskip .2in \text{and}  \hskip .2in s=\frac{x-y}{2}$$
then $$f(x,y)=F(x) + G(y) +k'(F(x)-G(y)).$$  If $x=y$ then
$f(x,0)=\gamma(x) =F(x) + G(x)$ and $F(x)=G(x)$.  We see that
$$f(x,y)=(1/2)(\gamma(x) + \gamma(y)) +(1/2) k'(\gamma(x) -
\gamma(y)).$$ 
Since we did not assume anything about how $f(x,y)$ was
constructed, only that it was split-complex holomorphic and
real-valued on a part of the real axis we can see that there is
only one extension.\end{proof}

\begin{example}
If $\gamma(t)=sinh(t)$, then using the power series expansion of $sinh(t)$, 
\newline $ f(t+k's) = sinh(t)cosh(s) + k'
cosh(t)sinh(s)$. Look at
\begin{eqnarray*}  f(x,y) &=& sinh\left(\frac{x+y}{2}\right)cosh\left(\frac{x-y}{2}\right) + k'cosh\left(\frac{x+y}{2}\right)sinh\left(\frac{x-y}{2}\right) \\
& = & (1/2)(sinh(x) + sinh(y)) + k'/2(sinh(x) - sinh(y)).\end{eqnarray*}
\end{example}

\vspace{0.3cm}

\begin{theo}\label{analyticcontinuation}
Let $f$ and $g$ be split-complex analytic functions on a domain $ D$ (open, connected
subset) of $\mathbb C'$.  If $f(z)=g(z)$ in a neighborhood of some $\alpha
\in D$ then $f=g$ in  $D$.\end{theo}

 We need two lemmas before we begin its proof.

\begin{lema}\label{coeff}
If $f(z)=\sum_{k=0}^\infty \gamma_k z^k$
then $\gamma_k = \frac{f^{(k)}(0)}{k!}, \ \ \ k=0,1,2....$
\end{lema}

This follows from \cite{Er} Lemma 1.5, which states that 
$f'(z) = \lim_{h\to 0} \frac{f(z+h) - f(z)}{h}.$  It is easy to see,
using this result that,  $f'(z) = \sum_{k=1}^\infty \gamma_k k  z^{k-1}$
and the formula follows by induction.

\begin{lema}\label{zerocoeff}
If $F(z)$ is split-complex analytic in a neighborhood of $\alpha\in D$ and
$F(z)\equiv 0$ in that neighborhood then $F^k(\alpha)=0$ for all $k$.
\end{lema}
Again we use the difference quotient definition to see that $$F'(\alpha) =
 \lim_{h\to 0} \frac{F(\alpha+h) - F(\alpha)}{h}=0$$ in the neighborhood. 
The higher derivatives are zero in the same way.

\vspace{.3cm}

{\bf Proof of Theorem \ref{analyticcontinuation}:}  This just follows the
one-variable proof given in Levinson and Redheffer's book \cite{LR}, 
page 147-8. In fact, set $F(z) = f(z) -g(z).$  Suppose there is a $\beta \in D$ with
$F(\beta)\ne 0.$  Join $\alpha$ to $\beta$ by a piecewise linear,
connected curve $z=\zeta(t)$ in $D$ with $a \le t \le b$, $\zeta(a) =
\alpha$ and $\zeta(b)=\beta.$  We define
$$S=\{t\in [a,b]\,|\, F^{(k)}[\zeta(t)]=0\,|\,k=0,1,2,...\}.$$
We see that $a\in S$ by Lemma \ref{zerocoeff} hence $S$ is not empty, so that $S$ has a least
upper bound $t_o$.  By definition we can take a sequence $\{t_j\} \subset S$
which converges to $t_o \in S$.  By continuity $\displaystyle
F^{(k)}(\zeta(t_0))=\lim_{j\to \infty} F^{(k)}(\zeta(t_j))=0,$ so that
$t_0 \in S.$  Using the series representation for $F$ in a neighborhood of
$\zeta(t_0)$ we see that $t_0$ cannot be the upper bound of $S$ and so no
$\beta $ exists.

\subsection{Timelike and spacelike Bj\"orling Problem for Lorentzian surfaces}

The classical Bj\"orling problem asks for the existence and uniqueness of a minimal surface in $\mathbb R^3$ that passes through a real analytic curve with a prescribed analytic unit normal along this curve.  Now, in this paper, since we are working with a Lorentzian metric we can study  two
forms of the Bj\"orling problem for minimal surfaces in $\mathbb
L^3$, namely, when the initial data $\gamma(t)$ is timelike or spacelike curve. In fact,   even though,  we can state a Bj\"orling problem when the initial data $\gamma(t)$ is a null curve, there is not uniqueness to the solution for this problem, as the following example shows.

\begin{example}\label{Exem}  
 Take the null cubic curve 
$$x(u)=(\frac{4}{3} u^3 + u, \frac{4}{3} u^3 - u, 2u^2),$$
with a unit normal field $N(u)=(2u,2u,1)$ along it.
Let 
\begin{eqnarray*}  y(v) &= \left(\frac{1}{2}sinh(2v), v,
\frac{1}{2}(cosh(2v) - 1) \right) \\
z(v) &= (v,v,0)\end{eqnarray*}
be two null curves. Then the two surfaces $f(u,v) = x(u) + y(v)$ and $g(u,v) = x(u) + z(v)$ are two
distinct  minimal surfaces containing the curve $x(u)$ with $N(u)$ a normal field
along  it.
\end{example}

Hence it remains only two problems to be studied.

\vspace{0.2cm}

Assume that $\gamma:I \to \mathbb L^3$ is  a regular analytic
timelike (respectively spacelike) curve in $\mathbb L^3$ and $W: I
\to \mathbb L^3$ is a unit analytic spacelike vector field
along $\gamma$ such that $\left <\gamma', W\right > =0$. The Bj\"orling problem is to determine a minimal
Lorentzian surface $X: \frak O \subset \mathbb C' \to \mathbb L^3$
such that $X(t,0) = \gamma (t)$ and $N(t,0) = W(t)$ (respectively
$X(0,s) = \gamma (s)$ and $N(0,s) = W(s)$ for all $s \in I$). In
our case, $\frak O $ is a split-complex domain with $I \subset
\frak O$ and $N :\frak O \to \mathbb L^3$ is the Gauss map of the
surface.

When $\gamma$ is timelike this problem is called the {\it timelike
Bj\"orling problem} and if $\gamma$ is spacelike, we call it the {\it
spacelike Bj\"orling problem}.

\vspace{0.2cm}

The following theorem describes the split-complex representation formula in the timelike Bj\"orling problem. We follow the notation established in Section 2 for $z = t + k's$ where $t$ and $s$ are conformal coordinates and note that the Lorentzian cross-product used in the theorem is defined by
$$\left<(u \times v),w\right> = det[u,v,w].$$ For its proof one follows the Proof of Theorem 3.1 in \cite{ACM}.

\begin{theo}(Timelike Bj\"orling representation)
Let $X:\frak O \subset \mathbb C'\to \mathbb L^3$ be a timelike
minimal surface in $\mathbb L^3$ and set $\gamma(t)=X(t,0)$, $W(t)
= N(t,0)$ on a real interval $I$ in $\frak O.$  Choose any simply
connected region $\frak R\subset \frak O$ containing $I$ over
which we can define split-holomorphic extensions $\gamma(z), W(z)$
for all $z\in \frak R$.  Then for all $z \in \frak R$ we have:
\begin{eqnarray}\label{eq:Bjorling}
X(z)=Re\left(\gamma(z) + k'\int_{t_o}^z W(w) \times \gamma'(w) dw
\right),
\end{eqnarray}
where $t_o$ is an arbitrary fixed point in $I$ and the path
integral is taken over any path in $\frak R$ from $t_o$ to $z$.
  \end{theo}
 \begin{proof}
 Since $X$ is a minimal immersion we can look at:
 $$\phi(z)=\frac{\partial X}{\partial z}(z),$$
 which is split-holomorphic over $\frak{O}.$
 We have $$X(z) =2 Re\int_{\gamma_z}\phi(w) dw$$
 with the constant of integration being the one that makes the expression $X(t,0)=\gamma(t)$ holds for all $t\in I.$
 We know that $N\times X_s = X_t$ so that:
 $$\phi(z)=\frac{1}{2}(X_t + k'X_s)=\frac{1}{2}(X_t + k'N\times X_t).$$
 From the definition of $\gamma, W$ we have $\phi(t,0)=\frac{1}{2}(\gamma'(t) + k'W(t)\times \gamma'(t))$.  This is a mapping from $I$ to $\mathbb C'^3$.  Then the argument in Proposition \ref{eq:3.2} shows that this has a  unique extension to:
 $$\phi(z) = \frac{1}{2}(\gamma'(z) + k'W(z)\times \gamma'(z)).$$  As in \cite{ACM} we end up with
 $$X(z)=Re\left( \gamma(z) + k'\int_{s_o}^z W(w) \times \gamma'(w) dw \right).$$ 
 \end{proof}

For the spacelike Bj\"orling problem the alternative choice of variable $w = k'z = s + k't$, described in the end of Section 2, is more convenient. It will allow us to get a spacelike Bj\"orling representation as follows.

\begin{theo} \label{eq:TB2} (Spacelike Bj\"orling representation)
Let $X:\frak O \subset \mathbb C'\to \mathbb L^3$ be a timelike minimal surface in $\mathbb L^3$ and set $\gamma(s)=X(0,s)$, $W(s) = N(0,s)$ on a real interval $I$ in $\frak O.$  Choose any simply connected region $\frak R \subset \frak O$ containing $I$ over which we can define split-holomorphic extensions $\gamma(w), W(w)$ for all $w\in \frak R$.  Then for all $w \in \frak R$ we have:
\begin{eqnarray}
X(w)=Re\left(\gamma(w) + k'\int_{s_o}^w W(\zeta) \times
\gamma'(\zeta) d\zeta \right),
\end{eqnarray}
where $s_o$ is an arbitrary fixed point in $I$ and the path
integral is taken over any path in $\frak R$ from $s_o$ to $w$.
\label{teospace}
\end{theo}

\begin{proof}
 Since $X$ is a minimal immersion we can look at:
 $$\psi(w)=\frac{\partial X}{\partial w}(w),$$
 is split-holomorphic over $\frak{O}.$
 We have $$X(w) =2 Re\int_{\gamma_w}\psi(\zeta) d\zeta$$
 with the constant of integration being the one that makes the expression $X(0,s)=\gamma(s)$ holds, for all $s \in I.$
 We know that $X_t=N\times X_s$ so that:
 $$\psi(w)=\frac{1}{2}(X_s + k'X_t)=\frac{1}{2}(X_s+ k'N\times X_s).$$
 From the definition of $\gamma, W$ we have $\psi (s,0)=\frac{1}{2}(\gamma' (s) + k'W(s)\times \gamma' (s))$.  This is a mapping from $I$ to $\mathbb C'^3$.  The argument in Proposition \ref{eq:3.2} shows that this has a  unique extension to:
 $$\psi(w) = \frac{1}{2}(\gamma' (w) + k'W(w)\times \gamma' (w)).$$  As in the previous case we end up with
 $$X(w)=Re\left( \gamma(w) + k'\int_{s_o}^w W(\zeta) \times \gamma' (\zeta) d\zeta \right).
 $$
 \end{proof}

 \begin{example} The Lorentzian helicoid of 3rd kind can be parametrized by $$X(t,s) = (sinh(t)cosh(s), sinh(t)sinh(s), s).$$
  Note that $sinh(t + k's)=sinh(t)cosh(s) + k' cosh(t)sinh(s)$ and 
  $cosh(t + k's) = cosh(t)cosh(s) + k'sinh(t)sinh(s).$
    Let $$\gamma(t)=(sinh(t),0,0),  \ \ \ \ W(t) = (0,-1/cosh(t), sinh(t)/cosh(t).$$  We also know that $(a+k'b)^{-1} = \frac{a-k'b}{a^2-b^2}$.
Thus
  $$\gamma(z) = (sinh(t)cosh(s) + k' cosh(t)sinh(s),0,0),$$
  $$\gamma'(z) = (cosh(t)cosh(s) + k'sinh(t)sinh(s),0,0),$$ and 
\begin{eqnarray*}W(z)= \left(0,  \frac{2(-cosh(t)cosh(s) + k'sinh(t)sinh(s))}{cosh(2s) + cosh(2t)}, \frac{sinh(2t)+ k' sinh(2s)}{cosh(2s) + cosh(2t)}\right).\end{eqnarray*}
Finally we see that $W(w) \times \gamma'(w) = (0, sinh(w),1),$ and 
$$Re\left( \gamma(z) + k' (0,cosh(z), z)\right) = (sinh(t)cosh(s), sinh(t)sinh(s), s).$$
\end{example}

\vspace{0.2cm}

Now let us show how to recover that Lorentzian helicoid through the spacelike Bj\"orling representation. 

In fact, 
 $\gamma(s)=X(0,s) = (0,0,s)$ and
$W(s)=N(0,s)=(- sinh(s),- cosh(s),0)$ are spacelike vectors. The extensions are 
$\gamma(w) = (0,0,w)$  and $ W(w)=(- sinh(w),- cosh(w),0).$ Then, we see that $W(\zeta) \times
\gamma' (\zeta) = (cosh(\zeta), sinh(\zeta),0).$  By using
Theorem \ref{teospace} we obtain:
\begin{eqnarray*}
X(s+k't)  &=&  Re\left( \gamma(s+k't) + k' \int_{s_0}^{s+k't} (cosh(\zeta), sinh(\zeta),0) \right) \\ &=& 
 (cosh(s) sinh(t), sinh(t) sinh(s),s).
\end{eqnarray*}

We observe that this Lorentzian helicoid is a ruled surface and it will be considered again in Example \ref{ex1}.

\vspace{0.2cm}

Since it is simple to move from the timelike solutions to the spacelike ones, in the following we will focus on the results corresponding to the timelike case.

The next result proves that the timelike Bj\"orling problem has a unique solution.

\vspace{0.2cm}

\begin{theo} \label{eq:theorem3}
There exists a unique solution to the timelike Bj\"orling problem
for minimal surfaces. In fact, if $\gamma, W$ are defined as in
the formulation of the timelike  Bj\"orling problem, then:

(1) there exists a simply connected open set $\frak O \subset
\mathbb C' $ containing $I$ for which $\gamma, W$ admit
split-holomorphic extensions $\gamma(z), W(z)$ over $ \frak O$ and
the mapping $X: \frak O \to \mathbb L^3$ given by
\begin{eqnarray}
X(z)=Re\left(\gamma(z) + k'\int_{t_o}^z W(w) \times \gamma'(w) dw
\right),
\end{eqnarray}
is a solution to the timelike Bj\"orling problem. Here $t_o \in
I$ is fixed but arbitrary.

(2)  If $X_1: \frak O_1 \subset \mathbb C'  \to \mathbb L^3,  \ X_2
:\frak O_2 \subset \mathbb C'  \to \mathbb L^3$, are two different
solutions to the timelike Bj\"orling problem, then $X_1$ and
$X_2$ coincide over the non-empty open set $\frak O_1 \cap \frak
O_2$.
  \end{theo}
 \begin{proof}

We start by proving (2).  The timelike
Bj\"orling representation shows that every solution of the
Bj\"orling problem is given by (\ref{eq:Bjorling}) on any simply
connected open set for which $\gamma(z)$ and $W(z)$ exist. So we can construct
the two split-holomorphic extensions, which are equal in a neighborhood of
$I$ in the plane $\mathbb C'$. It follows then from Theorem \ref{analyticcontinuation} that they agree in $\frak O_1 \cap \frak O_2.$

For (1),  let $\frak O \subset \mathbb C'$ be a open set such that
$I\subset \frak O$ and over which the split-holomorphic extensions
$\gamma(z), W(z)$ exist. We define the split-holomorphic mapping
$\phi: \frak O \subset \mathbb C' \to \mathbb C'^3:$
$$
\phi(z)=\frac{1}{2}(\gamma'(z) + k' W(z) \times \gamma'(z)).$$
 So, if $\phi = (\phi_1, \phi_2, \phi_3)$, it follows that
$$
- \phi_1(z)^2 + \phi_2(z)^2 + \phi_3(z)^2=0,$$ and $$ - |\phi_1(t,
0)|^2 + |\phi_2(t, 0)|^2 + |\phi_3(t, 0)|^2 = \frac{1}{4} ( 1 +
|W(t) \times \gamma'(t)|^2) > 0.
$$
Now we assume that $\frak O$ is simply connected and that
for all $z \in \frak O$, 
\newline $- |\phi_1(z)|^2 + |\phi_2(z)|^2 +
|\phi_3(z)|^2 > 0$.  Since $$2 Re \int_\gamma \phi_k dz =
\int_\gamma \left(\frac{\partial \psi_k}{\partial t}dt +
\frac{\partial \psi_k}{\partial s}ds \right) = \int_\gamma d\psi_k
=0,  $$ Theorem \ref{eq:immersion} assures us that
$$
 X(z) = 2 Re \int_{t_o}^z (\phi_1(w), \phi_2(w), \phi_3(w)) dw$$
 is a minimal immersion in $\mathbb L^3$, i.e,  $X: \frak O \subset \mathbb C' \to \mathbb L^3$ given by
 $$
X(z) = Re \left (\gamma(z) + k' \int_{t_{o}}^z W(w) \times \gamma'(w)
dw \right)
$$
is minimal surface.  Finally, $X$ satisfies the conditions of the
Bj\"orling problem. In fact, since $\gamma(z)$ and $W(z)$ are real
when restricted to $I$ , we have $X(t,0) = \gamma(t)$. Moreover,
one  has $$
 \frac{\partial X}{\partial t} (t,0) = \gamma'(t), \ \ \  \frac{\partial X}{\partial s} (t,0) = W(t) \times \gamma'(t),
 $$
which implies that
$$
W(t) \times \frac{\partial X}{\partial t} (t,0) = \frac{\partial
X}{\partial s} (t,0) = N(t, 0) \times \frac{\partial X}{\partial
t} (t,0) ,
$$
and so $N(t,0) = W(t)$.
\end{proof}

\begin{cor} \label{eq:corol1}
Let $\gamma: I \to \mathbb L^3$ be a regular analytic timelike curve in $\mathbb L^3$, and let $W:I \to \mathbb L^3$ be a spacelike analytic unit vector field along $\gamma$ such that $\left<\gamma', W\right > =0$. There exists a unique analytic minimal immersion in $\mathbb L^3$ whose image contains $\gamma(I)$ and such that its Gauss map along $\gamma$ is $W$.
\end{cor}

\begin{proof} It remains to prove the uniqueness since the existence comes from Theorem
\ref{eq:theorem3}.  Assume that  $X: M^2_1 \to \mathbb L^3$ is a
minimal immersion with a local isothermal coordinates system
$(U,\psi)$, where $U$ is an open set in $M^2_1$ and $\psi(U)=V$.  Choose $J \subset I$ so that $\gamma(J) \subset
X(U)$.  Locally $X|_U$ can be written as a minimal surface $\chi:
V \to \mathbb L^3$ defined by $X(\psi^{-1}(V)).$  There is an
$\alpha: J \to V$ such that $\chi(\alpha(t)) = \gamma(t)$ and
$N(\alpha(t))=W(t)$ for all $t \in J.$  We can see that $\alpha$
is  real analytic as follows.  The Jacobian of $\chi$ has rank 2.
At any point $\alpha(t_o) = p_o$ pick two coordinates such that
$(\chi_1, \chi_2)$ have invertible Jacobian at $p_o$.  Then
$\alpha(t) =(\chi_1, \chi_2)^{-1} \circ \gamma(t)$ is real
analytic and so has a split-holomorphic extension $\alpha(z): O
\subset \mathbb C' \to \mathbb C',$ where $O$ is open and $J
\subset O.$
 Writing $\alpha(t) = \alpha_1(t) + k' \alpha_2(t)$ and using the fact that $\gamma(t)$ is a 
regular curve we obtain $\alpha'^2_1 - \alpha'^2_2 \ne 0$.  Then one can
apply the inverse function theorem in a neighborhood of a point
$t_o \in J$ for which $\gamma(z)$ has non-null derivatives at $t_o$.
In fact, since the split-holomorphic complexification of the
real-analytic function $\alpha_j(t)$ is given by:
$$f_j(t,s) =\frac{1}{2}(\alpha_j(t+s) + \alpha_j(t-s)) + \frac{1}{2}k'((\alpha_j(t+s) - \alpha_j(t-s)),$$
the split-holomorphic extension is:
\begin{eqnarray*}f_1(t,s) + k'f_2(t,s) &= & 
\frac{1}{2}(\alpha_1(t+s) + \alpha_1(t-s)+\alpha_2(t+s) - \alpha_2(t-s))\\ & + & 
 \frac{k'}{2}((\alpha_1(t+s) - \alpha_1(t-s)+\alpha_2(t+s) + \alpha_2(t-s)),
\end{eqnarray*}
whose Jacobian determinant is $( \alpha_1'(t+s)+ \alpha_2'(t+s))(
\alpha_1'(t-s)-\alpha_2'(t-s))$. Hence, as the original curve
$\gamma(t)$ is timelike, for $s=0$ the determinant is non-zero.  Thus
we obtain a split-biholomorphic mapping $\alpha(z): A \subset
\mathbb C' \to B \subset \mathbb C'$, where $A$ is open subset of
$V$ which contains a real interval $(t_o - \epsilon, t_o +
\epsilon)$ and $B$ is an open subset of $V$. Hence the
minimal surface $X|_B: B \subset V \to \mathbb L^3$ can be
expressed as $\varphi: A \subset \mathbb C' \to \mathbb L^3$ with
$\varphi(z) = X(\alpha(z))$. Moreover, for all $t \in (t_o -
\epsilon, t_o + \epsilon)$ we have
\begin{eqnarray}
\varphi(t,0) = X(\alpha(t,0) = X(\alpha(t)) = \gamma(t),\\
N_\varphi(t,0)= N(\alpha(t,0) = N(\alpha(t)) = W(t).
\end{eqnarray}

Hence it follows from the uniquess of  $\varphi(z)$ that $X: M^2_1
\to \mathbb L^3$ is also unique. \end{proof}

Now let us consider the restricted timelike Bj\"orling problem: Let 
 $\gamma :I \to \mathbb L^3$ be a real analytic
curve in $\mathbb L^3$ with $<\gamma', \gamma'> = - 1$ and such
that $\gamma''(t)$ is spacelike for all $t \in I$. Construct a
minimal Lorentzian surface in $\mathbb L^3$ containing $\gamma$ as
a geodesic.

The next corollary, whose proof is similar to Corollary 3.5 in \cite{ACM}, gives the answer for the above problem.

\begin{cor} \label{eq:corol2} Let $\gamma: I \to \mathbb L^3$ be a constant speed analytic timelike curve in $\mathbb L^3$ such that $\gamma''(t)$ is spacelike for all $t \in I$. There exists a unique minimal Lorentzian immersion in $\mathbb L^3$ which contains $\gamma$ as a geodesic.
\end{cor}

Following \cite{ACM}, it is possible to construct examples of minimal immersions containg a given curve as geodesic. In the next example, we start with a  pseudo-circle in
$\mathbb L^3$, i.e, a planar timelike curve with non-zero constant
curvature.

\begin{example} \label{eq:pseudo} Any pseudo-circle contained in a timelike plane in $\mathbb L^3$ is congruent to a curve of the form $-x_1^2 + x_3^2 = R^2$, and may be parametrized by $\gamma(t) = R(sinh(t), 0, cosh(t))$.  It follows from Corollary \ref{eq:corol2} that there is a unique minimal immersion in $\mathbb L^3$ containing $\gamma$ as a geodesic. So, taking $W= - \frac{\gamma''}{|\gamma''|}$ i.e., $W(t) = - (sinh (t), 0, cosh (t))$ we get 
$$
\gamma(t) + k' \int^t W(\tau) \times \gamma'(\tau) d\tau = R(sinh(t), - k' t, cosh (t)).
$$
Hence the minimal immersion contaning $\gamma$ as geodesic, is
given by
$$
X(t,s) =  R(sinh(t)  cosh (s), - s, cosh (t)  cosh (s))$$ for $(t,s) \in
\mathbb R \times (- \frac{\pi}{2},   \frac{\pi}{2})$.
\end{example}

 We point out that this is a  surface of  revolution and it will be contained in formula (\ref{eq:rc}).

 Observe that if the minimal immersion in $\mathbb L^3$ contains a pseudo-circle
as a geodesic, the plane in which the pseudo-cirle is contained is
timelike. Hence we have a similar consequence to
Proposition 3.6 of \cite{ACM}, namely:

\begin{prop}\label{eq:prop1}
Any minimal timelike immersion in $\mathbb L^3$ containing a pseudo-circle
as a geodesic is congruent to a piece of a Lorentzian surface given by Example \ref{eq:pseudo}.
\end{prop}

For the spacelike Bj\"orling problem we obtain analogous results to Theorem \ref{eq:theorem3}, Corollaries \ref{eq:corol1}, \ref{eq:corol2} and Proposition \ref{eq:prop1}.

\section{Minimal timelike  surfaces of revolution}

Here we will give an alternative proof for the classification of timelike minimal surfaces of revolution in $\mathbb L^3$ given by Woestijne in \cite{W}, where one can also find the graphics of those surfaces. In our proof we show that those surfaces can be characterized as solutions of certain timelike or spacelike Bj\"orling problems.

We start by considering the different kinds of
 surfaces of revolution in $\mathbb L^3$, depending on the causal caracter 
of the axis of revolution, as obtained in \cite{Be}. They can be parametrized by:
\begin{equation}\label{eq:ra}
a) \ \ \ \ \ \ \ \   X(t,s) = (a(t),  b(t) \ cos(s), b(t) \ sin(s)),
\end{equation}
where 
$(a(t), b(t))$ is a timelike curve and $b(t)  \ne 0$.
\begin{equation}\label{eq:rb}
b) \ \ \ \ \ \ \ \   X(t,s) = (a(t) \ cosh (s), a(t) \ sinh (s), b(t)),
\end{equation}
 where $(a(t), b(t))$ is a timelike curve and $a(t)  \ne 0$.
\begin{equation}\label{eq:rc}
c) \ \ \ \ \ \ \ \  X(t,s) = (a(s) \ sinh(t), a(s) \ cosh (t), b(s)),
\end{equation}
 with $a(s) \ne 0$, $a'^2 + b'^2 \ne 0$.
\begin{equation}\label{eq:rd}
d) \ \ \ \ \ \ \ \ X(t,s) = (\frac{a(t) - b(t)}{\sqrt{2}} + \frac{a(t)
s^2}{2\sqrt{2}}  , \frac{a(t) + b(t)}{\sqrt{2}} -  \frac{a(t)
s^2}{2\sqrt{2}}, s a(t))
\end{equation}
  with $a'(t) b'(t) <0, \ a(t)  \ne 0$.

It is also known that all of these  surfaces  can be
conformally parametrized if the profile curves are parametrized
properly.
 
\vspace{0.3cm}
 Next we see examples of surfaces in $\mathbb L^3$ which will be necessary for the classification of timelike minimal  surfaces of revolution. We begin with  the following lemma.

\begin{lema} \label{eq:lema1}
let $\gamma(t)$ be a timelike analytic curve in  $\mathbb L^3$ contained in the timelike coordinate plane  $x_1, x_3$ or $x_1, x_2$-plane. Then  there exists a
unique timelike minimal immersion in $\mathbb L^3$, that
intersects orthogonally that plane along of $\gamma$, and is parametrized respectively by:
\begin{eqnarray}
a)  \ X(z)= (Re\; a(z), Im \int_{t_o} ^z \sqrt{a'^2 - b'^2} d\tau, Re\; b(z)), \ \text{if} \  \gamma(t) = (a(t), 0, b(t)),\\ 
b) \  X(z)= (Re\; a(z), Re \; b(z), Im \int_{t_o} ^z \sqrt{a'^2 - b'^2}
d\tau), \ \text{if} \   \gamma(t)= (a(t), b(t), 0).
\end{eqnarray}
  \end{lema}

 \begin{proof} In order to prove a) the Gauss map along the curve $\gamma$ is orthogonal
to $\gamma'$ and $e= (0,1,0)$. So, one has that $N(t) =
\frac{\gamma' (t) \times e}{| \gamma' (t) \times e|}$ and from
Corollary 3.1 we obtain the existence and uniqueness. Finally the
explicit formula above comes directly from the timelike Bj\"orling
representation. The proof of b) is similar.
 \end{proof}

\begin{example} ({\it Lorentzian elliptic catenoid}) Let $\gamma(t) = A (t, \ cos(t -\theta), 0)$ with $A>0$, $\theta \in \mathbb R$ and $t \in (\theta - \pi/2, \theta + \pi/2)$. Changing to a new parameter $u = t - \theta$ and 
using that $$sin(z) =  sin (t +k' s) = cos (s) \ sin (t) + k' \
cos (t) \ sin(s)$$ 
$$cos(z) = cos (t +k' s) = cos (t) \ cos (s) -  k'  \ sin (t) \ sin (s),$$ we have from Lemma \ref{eq:lema1}, a timelike minimal
surface which may be parametrized by
$$
X(u,v) = A(u + \theta, cos(u)  cos(v), cos(u)  sin (v)),  \ \ \text{where} \ \ (u,v) \in (-\pi/2, \pi/2) \times \mathbb R.
$$
\end{example}

\begin{example} ({\it Lorentzian hyperbolic catenoid}) Let $\gamma(t) = A (sinh(t+\theta), 0, t)$ with $A>0$, $\theta \in \mathbb R$ defined  for all $t > -\theta$.  Now,  setting $u = t + \theta$, one obtains from Lemma \ref{eq:lema1}  a timelike minimal
surface which may be parametrized by
$$
X(u,v) = A(sinh(u) cosh (v),  sinh (u) sinh (v),  u-\theta), \ \ \text{where} \ \ u>0, \ \ v \in
\mathbb R.
$$
\end{example}

\begin{example}( {\it Lorentzian surface with spacelike profile curve})
 Let 
 \newline $\gamma(s) = A (0,  cosh(s + \theta), s)$ where $A>0$, $\theta \in \mathbb R$ and $s > - \theta$. Choose a new parameter $v = s + \theta$ and consider a simple variation of Lemma \ref{eq:lema1},  for analytic spacelike curves  parametrized by $(0, a(s), b(s))$. Then taking $e= (1,0,0)$ and $N(s) = \frac{e \times \gamma'(s)}{| e \times \gamma'(s)|}$, we have the existence and uniqueness of the timelike minimal immersion given explicit by:
\begin{eqnarray}
X(w)= X(s+k' t) = ( Im \int_{s_o} ^w \sqrt{a'^2 + b'^2} d\tau, Re
\ a(w), Re \ b(w)) ,
\end{eqnarray}
where $w= k' z = s+ k' t$. By applying this result to the curve $\gamma(s) = A (0,  cosh(s + \theta), s)$ we get a
timelike minimal surface parametrized by
$$
X(u,v) = A (cosh(v)  sinh(u), cosh(v)  cosh(u), v - \theta), \ \ \text{where} \ \ v > 0, \ \ u
\in \mathbb R.
$$
\end{example}

The next example corresponds to the  Lorentzian parabolic catenoid.
To simplify the computations, we will use  a null frame of
$\mathbb L^3$ given by
$$
L_1= ( - \frac{\sqrt{2}}{2}, \frac{\sqrt{2}}{2}, 0), \ \ L_2 =  (
\frac{\sqrt{2}}{2}, \frac{\sqrt{2}}{2}, 0), \ \  L_3= ( 0, 0,
1).$$

\begin{example}( {\it Lorentzian parabolic catenoid})\label{eq:parab}
Applying Lemma \ref{eq:lema1}  to the analytic timelike curve $(p(t), q(t), 0)$
written with respect to the null frame $\{L_1, L_2, L_3\}$, we obtain existence and uniqueness of the timelike minimal immersion given by :
\begin{eqnarray}\label{eq:parabo}
X(z) = ( Re \ p(z), Re \ q(z), Im \int_{t_o} ^z \sqrt{-2 p'(\tau)
q'(\tau)} d\tau).
\end{eqnarray}
Applying this result to the curve 
$$\gamma(t) = A (\frac{1}{6} t^3 + \frac{B}{2} t^2 + \frac{B^2}{2} t, -(t+B), 0), \ \ \text{where} \ \ A>0, \ B \in \mathbb R, \  t < - B,$$
we get a timelike minimal surface
which may be parametrized by:
$$
X(t,s) = A (\frac{1}{6} t^3 + \frac{B}{2} t^2 + \frac{B^2}{2} t +
\frac{s^2}{2} ( t + B), -(t+B), s(t+B)),$$
with respect to the
null frame.
\end{example}

{\bf Proof of Theorem \ref{eq:rev}}: Consider a timelike minimal surface of revolution
parame-trized by the conformal immersion (\ref{eq:ra}):  In this
case, the $x_1, x_2$-plane intersects the surface orthogonally
along the curve $\gamma(t) = X(t,0) = (a(t), b(t), 0)$. Here
$N(t,0) \times \gamma' (t)$ is collinear with $e=(0,0,1)$,  and
from Bj\"orling representation one sees that the split-holomorphic
extensions $a(z)$ and $b(z)$ should satisfy
$$
Re\; a(z) = a(t),  \ \  \ \  Re\; b(z) = b(t) \;  cos (s).$$ Thus $\frac{\partial}{\partial t} Re\; a(z) =
\frac{\partial}{\partial s} Im\; a(z)$ and $\frac{\partial}{\partial s} Re\; a(z) =  \frac{\partial}{\partial
t} Im \; a(z).$ Since $Re\; a(z) = a(t)$, one obtains $a(t) = At
+ B$ for $A,B$ constants. Now applying the split-holomorphic
conditions for $b(z)$,  we find that $b(t) = C_1 \; cos (t) + C_2
\; sin (t)$, where $C_1, C_2$ are constants.  Since the immersion is
conformal, we have that $|\gamma'(t)|^2 = - b^2(t)$, which implies
that $A^2 = C_1^2 + C_2^2$ i,e., there is $\theta \in \mathbb R$
such that $C_1=  A \; cos (\theta)$ and $C_2 = A \; sin (\theta)$.
Substituting those values in $b(t)$ one has $b(t) = A cos (t -
\theta)$. Using Lemma 
\ref{eq:lema1}, the surface is a piece of the Lorentzian elliptic catenoid.
\vspace{0.2cm}

 Now, let us consider a minimal surface of revolution  parametrized by (\ref{eq:rb}), in which the rotation of the timelike analytic curve $(a(t), 0, b(t))$ is around the $x_3$-axis. Following a), we get $b(t) = At +B$  and  $a(t) = A  sinh (t + \theta)$, where $\theta \in \mathbb R$ and $A, B$ are constants. From Lemma \ref{eq:lema1} the resulting surface is congruent to a piece of the Lorentzian hyperbolic catenoid.

\vspace{0.2cm}
Now we will consider parametrization (\ref{eq:rc}). For that case, the $\gamma(s) = (0, a(s), b(s))$  is an analytic spacelike curve which is rotated around of the $x_3$-axis.
Following the same idea above, one obtains 
$$
Re \ a(w) = a(s) cosh(t), \ \ \ Re \ b(w) = b(s),
$$
where $w= k'z = s + k' t$.  Hence the split-holomorphic conditions
for both $a(w)$ and $b(w)$, imply that $b(s) = As + B$ and $a(s) =
C_1 cosh(s) + C_2 sinh(s)$, where $A,B, C_1, C_2$ are constants. As the immersion is conformal, $|\gamma'(s)|^2
= a^2(s)$, $C_1^2 - C_2^2 = A^2$, and so $a(s) = A \
cosh (s+\theta)$, $\theta \in \mathbb R$. Now  the surface obtained is
congruent to a piece of a Lorentzian surface with spacelike profile curve.

\vspace{0.2cm}

Finally  we consider a minimal surface of revolution parametrized
by the conformal immersion (\ref{eq:rd}) written as
$$
X(t,s) = (b(t) - \frac{a(t) s^2}{2},  a(t),  \ s a(t)),
$$
 with respect
to the null frame $\{L_1, L_2, L_3\}$. 
Then the $x_1, x_2$-plane intersects the surface orthogonally along
the curve $\gamma(t)=X(t,0) = (b(t), a(t), 0)$. We also obtain $N(t,0) \times
\gamma'(t)= (0, 0, \pm \sqrt{-2 a' b'})$ along $\gamma$. Hence
using representation (\ref{eq:parabo}) one gets that
$$
Re \ b(z) = b(t) - \frac{a(t) s^2}{2}, \ \ \ Re \ a(z) = a(t).
$$
Since $-2 a'(t) b'(t) = a^2$, it follows that $$a(t) = - A (t +B), \ \ \ 
b(t) = A ( \frac{1}{6} t^3 + \frac{B}{2}
t^2 + \frac{B^2}{2}t).$$
Applying the same variation of Lemma \ref{eq:lema1} used in Example 4.4, the surface is congruent to a piece of the Lorentzian parabolic catenoid.  $\square$
 
 \section{Minimal timelike ruled surfaces}

In this section we study the minimal timelike ruled surfaces in $\mathbb L^3$. Using our split-complex Bj\"orling representation, we will give an alternative proof to the classification obtained by Woestijne in \cite{W},  where the graphics of those surfaces can be found.

  \vspace{0.3cm} 
 We begin identifying the timelike ruled surfaces in $\mathbb L^3$, following Kim and Yoon in \cite{KY1} and \cite{KY2}.

 A ruled surface in $\mathbb L^3$ is defined by:
 \begin{equation} \label{eq:ruled}
 X(t,s) = \alpha(t) + s \beta(t), \ \ \ t\in J_1, \ s \in J_2,
 \end{equation}
with $J_1$ and $J_2$ open intervals in $\mathbb R$ and where $\alpha=\alpha(t)$ is a curve in $\mathbb L^3$ defined on $J_1$ and $\beta = \beta(t)$ is a transversal vector field along $\alpha$. The curve $\alpha = \alpha(t)$ is called the {\it base curve} and $\beta = \beta(t)$ the {\it director vector field}. In particular if $\beta$ is constant, the ruled surface is called cylindrical, and non-cylindrical otherwise.

 First, we suppose the base curve $\alpha$ is spacelike or timelike. In this case, the director vector field $\beta$ can be naturally chosen to be orthogonal to $\alpha$. In addition, since the ruled surface is timelike, we get different cases, depending on  the causal  character of the base curve $\alpha$ and the director vector field $\beta$, as follows:

{\bf Case 1} \ The base curve $\alpha$ is spacelike and $\beta$ is timelike.
In this case $\beta'$ must be spacelike since it lies in $[\beta]^\perp$. This
surface will be denoted by $X_+^3$.

 \vspace{0.2cm}

{\bf Case 2} \  $\alpha$ is timelike and $\beta'$ is non-null.
In this case the director vector field $\beta$ is always spacelike and the
surface will be denoted by $X_-^1$.

\vspace{0.2cm}

{\bf Case 3} \  $\alpha$ is timelike and $\beta'$ is lightlike.
In this case the director vector field $\beta$ is always spacelike and the 
surface will be denoted by $X_-^2$.

\vspace{0.2cm}

 But if the base curve $\alpha$ is a lightlike curve and the vector field $\beta$ along $\alpha$ is a lightlike vector field, then the ruled surface is called a {\it null scroll}. In particular, a null scroll with Cartan frame is said to be a $B$-scroll (\cite{KY1}, \cite{KY2}). It is also a timelike surface.
 
 \vspace{0.2cm}
We first give some examples of minimal timelike ruled surfaces.

 \begin{example} \label{ex1}{\it (Timelike helicoid of the 3rd kind)} Let
\begin{equation} \label{eq:X_+^3}
\begin{cases}
\gamma(t) = (t,0,0),\\
W(t) = \frac{1}{\sqrt{1 + c^2 t^2}}(0, -ct, 1).$$
\end{cases}
\end{equation}
Changing the parameter to $ct = sinh (u)$, one gets
$$
\gamma(u) = \frac{1}{c} (sinh (u), 0,0), \ \ \ \ W(u) =
\frac{1}{cosh (u)} (0, - sinh (u), 1).
$$
Using the timelike Bj\"orling representation, we obtain the solution
of the timelike Bj\"orling problem with respect to the given data $(\gamma, W)$, parametrized
by
$$
X(z) = \frac{1}{c} ( sinh(u) cosh (v), v, sinh (u) sinh (v)), \ \ \ z= u + k'v.
$$
\end{example}

\begin{example}\label{ex2} {\it (Timelike helicoid of the 1st kind)} Let us consider the data
\begin{equation} \label{eq:X_-^3}
\begin{cases}
\gamma(s) = (0,s ,0),\\
W(s) = \frac{1}{\sqrt{1 - c^2 s^2}}(cs, 0, 1).$$
\end{cases}
\end{equation}
Changing the parameter to $cs = sin (v)$, one gets
$$
\gamma(v) = \frac{1}{c} (0, sin (v), 0), \ \ \ \ W(v) =
\frac{1}{cos (v)} (sin (v), 0, 1).
$$
Applying the Bj\"orling representation to the spacelike curve $\gamma(s)$,
the solution of the Bj\"orling problem is
parametrized by:
$$
X(w) = \frac{1}{c} (u, sin(v) cos (u), sin (v) sin (u)), \ \ \ w= v + k'u.
$$
\end{example}

\begin{example}\label{ex3}{\it (Timelike helicoid of the 2nd kind)} Let
\begin{equation} \label{eq:X_-^32}
\begin{cases}
\gamma(s) = (0, s, 0),\\
W(s) = \frac{1}{\sqrt{c^2 s^2 - 1}}(1, 0, cs).$$
\end{cases}
\end{equation}
Changing the parameter to $cs = cosh(v)$, one gets
$$
\gamma(v) = \frac{1}{c} (0, cosh (v), 0), \ \ \ \ W(v) =
\frac{1}{sinh (v)} (1,0, cosh (v)).
$$
Hence the solution of the spacelike Bj\"orling problem is
parametrized by:
$$
X(w) = \frac{1}{c} (cosh (v) sinh (u), cosh (v) cosh (u), u), \ \ w = v + k'u.
$$
\end{example}

\begin{example}\label{ex4} {\it (Conjugate of Enneper`s timelike surface of the 2nd kind)} Let \begin{equation} \label{eq:X_-^323}
\begin{cases}
\gamma(s) = (0, s, 0),\\
W(s) =  \frac{1}{\sqrt{1 - 2 sc}}(cs, 0, 1-cs)$$
\end{cases}
\end{equation}
with $c \ne 0$ and $1 - 2sc > 0$. Changing the parameter to $v^2 = 1- 2 cs$, $v >0$, the solution of the  Bj\"orling
problem is parametrized by:
$$
X(w) = - \frac{1}{6c} (3u + 3u v^2 + u^3, \ 3v^2 + 3 u^2 -3, \ 3u
- 3u v^2 - u^3), \ \ \ w= v + k'u.
$$
\end{example}

\begin{example}\label{ex5} ({\bf $B$-scroll})
Let $\alpha = \alpha (t)$ be a lightlike curve in $\mathbb
L^{3}$ with Cartan frame $\{A,B,C\}$ i.e., $A, B, C$ are vector
fields along $\alpha$ in $\mathbb L^{3}$ satisfying the
following conditions:
\begin{equation} \begin{array}{rl}\label{cond1}
&\left< A,A\right> = \left< B,B\right> =0,\, \left< A,B\right> =1,\\
&\left< A,C\right> = \left< B,C\right> =0,\,\left< C,C\right>=1, \\
&\alpha'= A, \ \ C'= -a A - c(t) B,
\end{array}
\end{equation}
where a is a constant and $c(t)$ a nowhere vanishing function.

The surface defined by  $X(t,s) = \alpha(t) + s B(t)$ is a timelike surface in
$\mathbb L^3$ called a $B$-{\it scroll}. Following  \cite{W}  a $B$-scroll is minimal if and only if it is flat, i. e., $B'(t) \equiv 0$ and $C'= - c(t) B$. 

It is possible to study this surface in the context of  timelike (spacelike)  Björling problem. In fact, let us reparametrize it by taking   the curve $\gamma (t)= \alpha(t) + s(t) B(t)$ with $s'(t)<0$ ($s'(t)>0$).  Then $\left< \gamma'(t),\gamma'(t)\right> = s'(t)<0 (>0)$ and $\gamma (t)$ is a timelike (spacelike) curve. In order to simplify the computations, take  $s(t)=-t \, (s(t)=t)$. Now $\gamma (t)=  \alpha(t)-tB(t)$ ($\gamma (t)=  \alpha(t)+tB(t))$  and $W(t)=C(t)$ are the timelike (spacelike) Björling data. Using (\ref{cond1}) we have $\left<\gamma'(t),W(t)\right>=0$. Using formula (\ref{eq:Bjorling}) we obtain the parametrization  of the  timelike (spacelike) Bj\"orling problem. 

For instance, taking the lightlike curve $\alpha(t)=\displaystyle(\frac{-t^3}{6 \sqrt 2}-\frac{t}{\sqrt 2},\frac{-t^2}{2},\frac{-t^3}{6 \sqrt 2}+\frac{t}{\sqrt 2})$   and the lightlike vector field $B(t)=\displaystyle(\frac{1}{\sqrt 2},0,\frac{1}{\sqrt 2})$ we obtain  the $B$-scroll
$$X(s,t)= \alpha(t)+s B(t)=\displaystyle(\frac{-t^3}{6 \sqrt 2}-\frac{t}{\sqrt 2}+ \frac{s}{\sqrt 2},\frac{-t^2}{2},\frac{-t^3}{6 \sqrt 2}+\frac{t}{\sqrt 2}+ \frac{s}{\sqrt 2}).$$
\end{example}

Using the reparametrization given above, we obtain the Björling data:
$$\gamma(t)=\left(\frac{-t^3}{6 \sqrt 2}-{\sqrt 2}\,t,\frac{-t^2}{2 },\frac{-t^3}{6 \sqrt 2}\right),$$
and $$W(t)=C(t)=A(t) \times B(t)=(\frac{t}{\sqrt 2},1,\frac{t}{\sqrt 2}).$$
After  once more using formula (\ref{eq:Bjorling}), we obtain the surface  parametrized by $X(t,s)=(X_1(t,s),X_2(t,s),X_3(t,s))$, where

\begin{equation}\begin{array}{rl}
&X_1(t,s)=\displaystyle -\frac{s^3+3 s^2 t+3 s t^2+t^3+12 t}{6
   \sqrt{2}}; \\
&X_2(t,s)= \displaystyle -\frac{ (s+t)^2}{2}; \\\\
&X_3(t,s)=\displaystyle -\frac{s^3+3 s^2 t+3
   s \left(t^2-4\right)+t^3}{6 \sqrt{2}}.
\end{array}
\end{equation}
Observe that setting $s=0$ in the parametrization above we get the curve $\gamma(t)$, as expected.

\vspace{.3cm}
 
{\bf Proof of Theorem \ref{eq:ruled3}} We follow closely the proof of Theorem 6.1 in \cite{ACM} and consider all the possible cases, depending on the causal character of the base curve and director vector field.  

\vspace{0.2cm}
 
 {\bf Case 1} Let $X$ be a non-cylindrical ruled surface of type $X_+^3$, parametrized by (\ref{eq:ruled}) such that $\left<\beta, \beta\right> = -1$ and $\left<\alpha', \beta'\right> = 0$. In
this case $\alpha$ is the striction curve and the parameter is the arc-length on the curve $\beta$. We define the
distribution parameter as
\begin{equation} \label{eq:para}
\lambda(t) = \frac{\left<\alpha' \times \beta, \beta'\right >}{\left<\beta',
\beta'\right>},
\end{equation}
since $\left <\alpha' \times \beta, \beta'\right> \ne 0$. In fact, 
$\alpha'\times \beta = \lambda \beta'$ and $X_t \times X_s =
\lambda \beta'+ s \beta'\times \beta$. Moreover $|| X_t \times
X_s||^2 = (\lambda^2 + s^2) \left <\beta', \beta'\right>$. So the striction
curve is a curve on the surface, obtained by setting $s=0$. The
Gauss map on the ruled surface is $$N(t,s) = \frac{\lambda(t)
\beta'(t) + s \beta'(t) \times \beta(t)}{\sqrt{\lambda^2(t) + s^2}
\ ||\beta'(t)||}.$$ 
Hence $\left<N(t_o, 0), N(t_o, s)\right> = \frac{1}{\sqrt{1 + c^2 s^2}}$, where $t_o
\in J_1$ and $c = \frac{1}{\lambda(t_o)}$. Let us assume $N(t_o,0)
= (0,0,1)$ and $L_s:= X(t_o, s)$, $s\in J_2$ parametrizing
the $x_1$-axis. Since $\left<N,N\right>=1$ we can assume $$N(t_o,s) =
\frac{(0, -cs, 1)}{\sqrt{1+ c^2 s^2}}$$ and that $L_s$ is
parametrized as $\gamma(s) = (s,0,0)$. So the minimal timelike 
ruled surface is the solution to the Bj\"orling problem with the
data: $\gamma(s) = (s,0,0), W(s) = \frac{1}{\sqrt{1 + c^2
s^2}}(0, -cs, 1).$ Hence the surface is
a piece of the timelike helicoid of the 3rd kind, according to Example \ref{ex1}.

{\bf Case 2} In this case the vector $\beta'$ is assumed to be 
non-null, so we must consider two subcases depending on whether
$\beta'$ is spacelike or timelike. In any subcase we have the
parametrized surfaces given by (\ref{eq:ruled}), with $\left<\beta, \beta\right> = 1$ and $\left <\alpha', \beta'\right> = 0$ and $\alpha$ is the striction curve. We define the
distribution parameter by (\ref{eq:para}) and conclude that $||X_t
\times X_s||^2 = (\lambda^2 - s^2) \left <\beta', \beta'\right>$. It follows
that if $\beta'$ is spacelike the striction curve
$\alpha$ is on the surface. If $\beta'$ is
timelike, we cannot have $s=0$.

 a) If $\beta'$ is spacelike, one gets that
 $$N(t,s) = \frac{\lambda(t) \beta'(t) + s \beta'(t) \times \beta(t)}{\sqrt{\lambda^2(t) - s^2} \ ||\beta'(t)||}.$$ Hence $\left <N(t_o, 0), N(t_o, s)\right> = \frac{1}{\sqrt{1 - c^2 s^2}}$ where $t_o
\in J_1$ and $c = \frac{1}{\lambda(t_o)}$. Now we assume $N(t_o,0)
= (0,0,1)$ and that $L_s:= X(t_o, s)$ for $s \in J_2$ parametrizes
the $x_2$-axis. Since $\left <N,N\right>=1$, we can assume $$N(t_o,s) =
\frac{(cs, 0, 1)}{\sqrt{1 - c^2 s^2}}$$ and $L_s$ is
parametrized by $\gamma(s) = (0,s,0)$. Hence the timelike minimal
ruled surface is the solution to the Bj\"orling problem with the
data: $\gamma(s) = (0,s,0), W(s) = \frac{1}{\sqrt{1 - c^2
s^2}}(cs, 0, 1).$ According to Example \ref{ex2}, the ruled surface
is a piece of the timelike helicoid of the 1st kind.

\vspace{0.3cm}

  b) If $\beta'$ is timelike, one gets that
 $$N(t,s) = \frac{\lambda(t) \beta'(t) + s \beta'(t) \times \beta(t)}{\sqrt{s^2 - \lambda^2(t)} \sqrt{|\left <\beta'(t), \beta'(t)\right>|}}.$$ Hence $\left <N(t_o, s_o), N(t_o, s)\right> = \frac{- \lambda^2(t_o) + s_o
s}{\sqrt{s_o^2 - \lambda^2(t_o)} \sqrt{s^2 - \lambda^2(t_o)}}$,
where $t_o \in J_1$ and $s_o \in J_2$ fixed. Now we assume $N(t_o,
s_o) = (0,0,1)$ and $L_s:= X(t_o, s)$, $s \in J_2$,
parametrizes the $x_2$-axis. Since $\left <N,N\right>=1$ we can assume
$$N(t_o,s) = \frac{1}{\sqrt{s_o^2 - \lambda^2(t_o)} \sqrt{s^2 -
\lambda^2(t_o)}} (|\lambda(t_o)| (s_o -s), 0, - \lambda^2(t_o) + s_o
s)$$
 and $L_s$ is parametrized by $\gamma(s) = (0,s,0)$.

 Now we take $s_o = \sqrt{2} |\lambda|$ and, substituting in $N(t_o, s)$, one gets
 $$
 N(t_o, s) = \frac{1}{\sqrt{s^2 - \lambda^2}}(- s + \sqrt{2} |\lambda|, 0, - |\lambda| + \sqrt{2} s),
 $$
where $\lambda$ is calculated at $t_o$. Now
composing with the orthogonal transformation of $\mathbb L^3$,
written in the canonical coordinates as:
$$
\begin{pmatrix}
-\sqrt{2} & 0 &-1\\
0& 1& 0\\
-1& 0 & -\sqrt{2}
\end{pmatrix}
$$
we obtain $$N(t_o,s) = - \frac{(1, 0, cs)}{\sqrt{c^2 s^2 -1}},
$$
where $c = \frac{1}{|\lambda|}$. Hence the timelike ruled surface
is the solution of the Bj\"orling problem with respect to the data: $\gamma(s) =
(0,s,0), W(s) = - \frac{1}{\sqrt{c^2 s^2-1}}(1, 0, cs)$. According
to Example \ref{ex3}, it is a piece of the timelike
helicoid of the 2nd kind.

\vspace{0.3cm}

 {\bf Case 3} Let $X$ be a non-cylindric ruled surface of type $X_-^2$, which may be parametrized by (\ref{eq:ruled}) where $\left <\alpha', \alpha'\right> =-1$, $\left <\alpha', \beta\right> =0$,
$\left <\beta, \beta\right> =1$ and $\left <\beta', \beta'\right> =0, \ (\beta'\ne 0)$.

 Consider the non-zero smooth functions $$-||X_t||^2(t_o,s) = 1 - 2s \left <\alpha'(t_o), \beta'(t_o)\right>$$
 and $\left <\beta'\times \beta, \alpha'\times \beta\right>(t_o) = - \left <\beta', \alpha'\right >(t_o)$. As $\beta \times \beta'= \beta'$, 
we have $$N(t_o, s) = \frac{1}{\sqrt{1- 2sc}} (\alpha' \times \beta -
s \beta')(t_o),$$ where $c = \left <\alpha'(t_o), \beta'(t_o)\right >$. Moreover $\left <N(t_o, s), N(t_o, 0)\right > = \frac{1- sc}{\sqrt{1-2sc}}.$
So, one may assume that $N(t_o, 0) = (0, 0, 1)$ and $X(t_o, s)$, $s \in J_2$, parametrizes the $x_2$-axis. Since 
$\left <N,N\right >=1$, it follows that
$$N(t_o, s) = \frac{1}{\sqrt{1- 2sc}} ( cs, 0, 1 - sc).
$$
So this timelike ruled surface is a solution of the spacelike Bj\"orling
problem with respect to the data: $\gamma(s) = (0,s,0), W(s) =
\frac{1}{\sqrt{1 - 2 sc}}(cs, 0, 1-cs).$ It corresponds to a piece of
the conjugate of Enneper`s timelike surface of the 2nd kind, just
as in Example \ref{ex4}.

\vspace{0.3cm}

{\bf B-scrolls}
Each of the cases above has essentially one surface, but the class of B-scrolls is larger, so we must use a different proof, which is similar to the proof found in \cite{W}.  We will find a simple representation for the B-scroll using the Bj\"orling procedure.

Begin with a timelike ruled surface $f(t,s) = \alpha(s) + t \beta(s)$, where $\left< \alpha'(s), \alpha'(s) \right> =0$ and $\left< \beta(s),\beta(s) \right> = 0$. This gives:
\begin{eqnarray*}
f_s&=&\alpha' + t\beta'\\
f_t&=&\beta.
\end{eqnarray*}
Since the surface is timelike we must have $\left< \alpha'(s), \beta(s) \right> \ne 0.$  We can first find a multiple of $\beta(s)$ so that  $\left< \alpha'(s), \beta(s) \right> =1.$
We construct a pseudo-orthonormal frame along $\alpha(s)$ using $\{\alpha',  \beta, m=\alpha' \times \beta\}$.  From the inner products we see there are functions $\{x_1(s), x_2(s), x_3(s)\}$ along the curve $\alpha(s)$ so that 
\begin{eqnarray}
\alpha'' &=& x_1 \alpha' + x_3 m\\
\beta' &=& -x_1 \beta + y_3 m\\
m' &=& -y_3 \alpha' -x_3\beta.
\end{eqnarray}
The surface unit normal is $f_s\times f_t=m+t \beta'\times \beta$.  We note that 
$\beta'\times \beta = y_3m \times \beta$ is a multiple of $\beta$, say $d(s)\beta$.  Thus the surface normal  is $N(s,t) = m + t d(s) \beta.$  $N_s = m' + t d' \beta + td \beta' = -y_3 \alpha' -x_3\beta+ t d' \beta + td \beta'  $ and must be a linear combination of $f_s$ and $f_t$.   Thus $d=-y_3$.  The shape operator has the form:
$\begin{pmatrix} y_3 & 0\\ *& y_3
\end{pmatrix}$, so by minimality $y_3=0$.  Finally we see that $\alpha$ is a pre-geodesic, and, by reparametrizing the curve we get $x_1=0$.  Thus our surface is a B-scroll as in Example 5.5.  We can find its (simple) spacelike Bj\"orling represention using $\gamma(s) =\alpha(s) + s \beta$ and $W(s) = m(s)$.  $W \times \gamma' =( \alpha' \times \beta) \times (\alpha' + \beta)$ = $(\alpha'-\beta)$.
\begin{eqnarray*}X(w) &=& Re\left( \gamma(w) + k' \int_{s_o}^w (\alpha'-\beta) d\zeta \right )\\&=&Re(\alpha(w))+ Im(\alpha(w))+ Re(\beta)(s-t) +\,Im(\beta)(t-s) .  \  \square  \end 
{eqnarray*} 

\vspace{0.2cm}
{\bf Acknowledgements} 
 The second author's research has been supported by a CAPES Grant/Brazil. She thanks  the Mathematics Department at University of California at Irvine for its hospitality. The third author would like to express his thanks to the Institute of Mathematics and Statistics (IME) at the University of  S\~ao Paulo S.P. for its hospitality and the State of S\~ao Paulo Research Foundation (FAPESP) for financial support during the development of part of this work.

\end{document}